\documentclass[review]{elsarticle}
\usepackage{amssymb}

\usepackage{thumbpdf} 

\newtheorem{Theorem}{Theorem}

\newtheorem{Corollary}[Theorem]{Corollary} 
\newtheorem{Proposition}[Theorem]{Proposition}

\newtheorem{ex}[Theorem]{Example}
\newtheorem{remark}[Theorem]{Remark}

\newcommand{\bR}{{\mathbb{R}}}

\title{Intersection bodies that are not polar zonoids: A flat top condition in dimensions four and six}
\author[M. A. Alfonseca]{M. Angeles Alfonseca}
\address{Department of Mathematics, North Dakota State University}
\ead{maria.alfonseca@ndsu.edu}

\tnotetext[t1]{Work partially supported by the NSF grant DMS-1100657.} 
\tnotetext[t2]{AMS 2010  Subject Classification: 52A21, 52A20. Keywords: Intersection bodies, zonoids, bodies of revolution}

\begin{document}

\begin{abstract} 
We find general geometric conditions on a convex body of revolution $K$, in dimensions four and six, so that its intersection body $I\!K$ is not a polar zonoid. We exhibit several examples of intersection bodies which are are not polar zonoids. 
\end{abstract}
\maketitle

\section{Introduction}

Let $K,L$ be origin-symmetric star bodies in $\bR^n$. The body $K$ is called the {\it intersection body of $L$}, and denoted by  $K=I \! L$, if for every $\xi \in S^{n-1}$, the radial function of $K$ in the direction $\xi$ is equal to the $(n-1)$-dimensional volume of the central section of $L$ perpendicular to $\xi$, {\it i.e.}  $ \rho_K(\xi)={\mbox{Vo}}l_{n-1}(L\cap \xi^\perp)$. By using the  spherical Radon transform $R$  (see \cite{G}, page 429), the above relation can be written as
\[
      \rho_K(\xi)=\frac{1}{n-1} \int_{S^{n-1}} \rho_L^{n-1} (\theta) \, d \theta= \frac{1}{n-1}\, R( \rho_L^{n-1})(\xi),           
      \qquad \forall \xi\in S^{n-1}.
\]
{\it Intersection bodies of star bodies} were introduced by Lutwak in \cite{Lu}. A more general class of {\it intersection bodies} is defined as the closure in the radial metric of the class of intersection bodies of star bodies. Intersection bodies proved to be crucial to solve the long-open Busemann-Petty problem (\cite{GKS}, see also \cite{G,K} for the history of the solution). The class of intersection bodies contains the class of polar zonoids, as proved by Koldobsky in \cite{K2} (see also \cite{KY} Chapter 6), and this containment is strict \cite{K1,KK1}. Recall that a compact convex set $K\subset \mathbb{R}^n$ is a zonoid if its support function is the cosine transform of a finite even Borel measure on $S^{n-1}$, {\it i.e.} for all $u\in S^{n-1}$
\[
       h_K(u)=C(\mu) (u):=\int_{S^{n-1}} |u \cdot v| d\mu(v),
\]
where $C$ denotes the cosine transform. More information about zonoids and their properties can be found, for example, in the paper by Bolker \cite{B}, or in Chapter 4 in  Gardner's book \cite{G}.

In the paper \cite{Sch3},  Schneider and Wieacker conjectured that for most convex bodies $K$, the intersection body of $K$ is not a polar zonoid (see also \cite{Sch4}). Here, ``most" is understood in the Baire category sense. Not much progress has been achieved on this conjecture so far, although there exist many examples of intersection bodies that are not polars of zonoids, such as the unit balls of certain Banach subspaces of $L_p(\mathbb{R}^n)$, for $0<p<1$ and $n \geq 3$ (see  \cite{K1,KK1,KZ}). One of the goals of the present work is to provide geometric criteria to construct examples of intersection bodies of revolution that are not polar zonoids. 

\begin{figure}[h!]
 \begin{center}
 \includegraphics[width=4.1in]{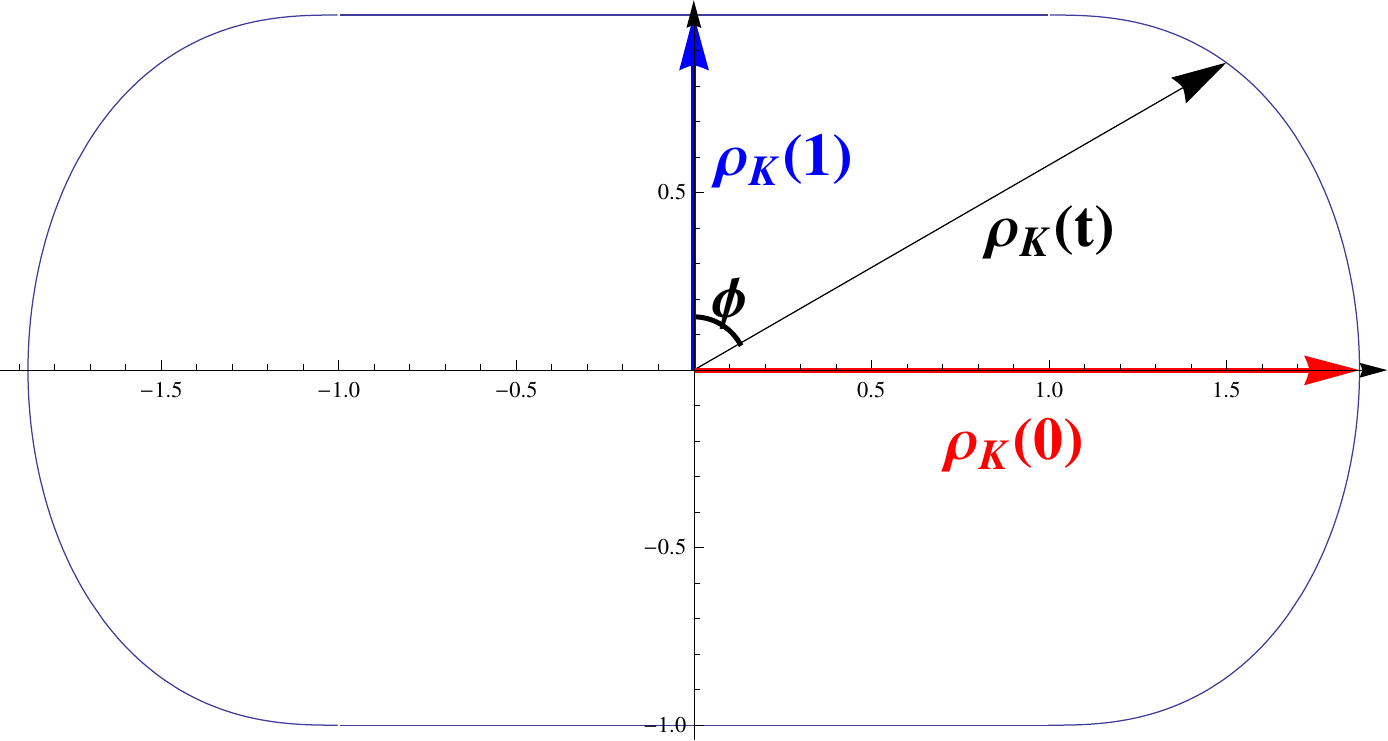}
 \end{center}
\caption{\small Definition of $\rho_K(t)$, where $t=\cos\phi$. }
\label{raf}
\end{figure}

The main result of the paper is a general geometric condition on a convex body of revolution $K$ that guarantees that $I\!K$ is not a polar zonoid, in dimensions four and six. As far as the author is aware, this is the first such general geometric condition in the literature. Taking into account that the polar of a zonoid is an intersection body, we follow the approach used by Lonke in \cite{L}, where he constructed a four-dimensional zonoid with a face, whose polar is a zonoid. If $K \subset \mathbb{R}^n$  is an origin-symmetric convex body of revolution around the $x_n$ axis, then its radial function $\rho_K$ is  {\it rotationally symmetric}, {\it i.e.} it can be defined as a function of $t$, the cosine of the vertical angle in spherical coordinates, by  $\rho_K(\sqrt{1-t^2}\,\xi,t)=\rho_K(t)$, $0\leq t \leq 1$, for all $\xi \in S^{n-2}$. Figure \ref{raf} shows a two dimensional section of a body $K$, where the axis of revolution is the vertical axis, and the radius of $K$ in the direction that forms an angle $\phi$ with the axis of revolution is $\rho_K(t)$, where $t=\cos \phi$. Because of the origin and rotational symmetries, knowing $\rho_K(t)$ for $0 \leq t \leq 1$ completely determines the body $K$. The value $t=1$ corresponds to the axis of revolution, and the value $t=0$ corresponds to the {\it equator} of the body. 

\bigskip

 With this notation, our four-dimesional result is the following.\\

 {\bf Proposition 1:}
{\it If $K$ is a four-dimensional origin-symmetric convex body of revolution whose radial function $\rho_K(t)$ is of class $C^2$ in a neighborhood of $t=1$, and  }
\begin{equation}
   \label{introprop}
         2\rho_K^4(1)>3 \left( \int_0^1  \rho_K^{3}(t)  \, dt  \right) \left( \rho_K(1)+\rho_K'(1) \right),
\end{equation}
{\it then the intersection body of $K$ is not a polar zonoid.} \\

The equation in Proposition 1 has an important geometric meaning. It reflects the fact that the body $K$ is ``flat" at the axis of revolution, in the following sense: Let us change our point of view for a moment, and consider the non-negative, even function $f_K(x)$,  $x\in [-\rho_K(0),\rho_K(0)]$ that defines the boundary of $K$, {\it i.e.} that every point on the boundary can be written either as $(x,f_K(x))$ or $(x,-f_K(x))$. If we assume that $\rho_K$ has a continuous second derivative everywhere, then so does the function $f_K(x)$, and  $f_K''(0)=-\frac{\rho_K(1)+\rho_K'(1)}{\rho_K(1)^2}$. Thus, Proposition 1 gives a bound for the second derivative of $f_K$ at $x=0$ (at the axis of revolution of $K$). In particular, the condition trivially holds if $\rho_K(1)+\rho_K'(1)=0$  (Corollary  \ref{flattop4}). We will refer to a body $K$ satisfying $f_K''(0)=0$, or  $\rho_K(1)+\rho_K'(1)=0$, as a {\it  body with a flat top}, or a body that is {\it flat at the axis of revolution}.

Note that the statement of Proposition 1 asks for $\rho_K$ to be of class $C^2$, but then condition (\ref{introprop}) involves only its first derivative. This is due to our choice of the variable $t$ for the radial function. If we let $\widetilde{\rho_K}(z)=\rho_K(t)$, where $z=\sin \phi$ and $t=\cos \phi$, then the term $\rho_K(1)+\rho_K'(1)$ corresponds to $\widetilde{\rho_K}(0)+\widetilde{\rho_K}''(0)$, and we see that the $C^2$ hypothesis around the axis of revolution is indeed needed. The proof of Proposition 1 will show this more clearly.

 In dimension six, the result is more involved:\\

{\bf Proposition 4:}
{\it Let $K$ be a six-dimensional origin-symmetric convex body of revolution with radial function $\rho
_K(t) $ of class $C^2$ in a neighborhood of  $t=1$.  Let $\displaystyle h(x)= \int_0^x \rho_K^{5}(t) (x^2-t^2)  \, dt$,  $\displaystyle k(1)=\int_0^1 t^2 \rho_K^5(t) \,dt$, and  $r(t)=\rho_K^5(t)$. If 
\begin{equation}
 \label{six}
   h^2(1)\left(5r(1)+r'(1) \right)+ 24k^3(1) < 12h(1)k(1)r(1),
\end{equation}
then $I\!K$ is not a polar zonoid.}\\

Proposition 4 shows that the situation in dimension six is different than in dimension four. In the six-dimensional case, having a flat top  (which means that the term  $((5r(1)+r'(1))$ in equation (\ref{six}) is equal to zero)  is not enough to guarantee that $I\!K$ is not a polar zonoid, and we need $K$ to satisfy also the condition $2k^2(1)<h(1)r(1)$. The geometric meaning of this condition and several applications of it are the subject of Section 4.

 It must be noted that our problem is more interesting starting from dimension five, for the following reason: In dimensions three and four, there exist convex bodies that are not zonoids (such as the octahedron), while every origin-symmetric convex body is an intersection body, as proved by Gardner in the three dimensional case  \cite{G2}, and by Zhang in dimension four \cite{Zh}.  However, there are two reasons why the four dimensional case is worth studying. First, while all four-dimensional convex bodies are intersection bodies, they are not necessarily {\it intersection bodies of  convex bodies}, which is the class of bodies we are considering in this paper. Secondly, the formulas of the Radon and cosine transform of rotationally symmetric functions are especially simple in dimension four, which allows us to identify the geometric condition in this case (Section 3), and then apply the same method in dimension six (Section 4). As expected, the conditions are more involved in dimension six than in dimension four. In Appendices A-C we specifically compute several examples, as applications of Proposition \ref{flattop4full} and Corollary \ref{flattop6}. Even in the four dimensional case the computations are  hard,  and were done using {\it Mathematica}.

\section{The spherical Radon transform and the cosine transform for rotationally symmetric functions}

Let $K \subset \mathbb{R}^n$  be an origin-symmetric convex body of revolution around the $x_n$ axis.
Its radial function $\rho_K$ is  {\it rotationally symmetric}, {\it i.e.} it can be defined as a function of $t$, the cosine of the vertical angle in spherical coordinates, by $\rho_K(\sqrt{1-t^2}\,\xi,t)=\rho_K(t)$, $0\leq t \leq 1$, for all $\xi \in S^{n-2}$. The intersection body of $K$ is defined by the relation $\displaystyle \rho_{I\!K}=\frac{1}{n-1}R(\rho_K^{n-1})$. By Busemann's theorem \cite{Bu},  $I\!K$ is convex. By the rotationally invariance of $R$, $I\!K$ is a body of revolution (see \cite{G}, Appendix C.2). 

Given a convex body $Z \subset \mathbb{R}^n$, its support function  $h_Z$  is defined  by $h_Z(x)=\sup \{ \langle x,u \rangle: u\in Z \}$, for $x\in \mathbb{R}^n$. It was shown by Bolker \cite{B} that a convex body $Z$ is a zonoid centered at the origin if and only if its support function can be represented in the form 
\begin{equation}
 \label{cosine}
     h_Z(x)=\int_{S^{n-1}} |\langle x,v \rangle| d \mu(v),
\end{equation}
where $\mu$ is a non-negative even measure. The integral operator on the right hand side of (\ref{cosine}) is called the cosine transform, which we will denote by $C$, thus  rewriting  (\ref{cosine}) as $h_Z=C\mu$. More information about zonoids and the cosine transform can be found in \cite{Sch} Section 3.5, and in \cite{G}, Chapter 4 and Appendix C.2. 

Our goal is to find conditions on a convex body $K$ to determine whether its intersection body $I\!K$ is or is not a polar zonoid. Thus, we need to study whether or not there is a non-negative even measure $\mu$ on $S^{n-1}$, such that $h_{(I\!K)^*}=C \mu$, where $(I\!K)^*$ denotes the polar body of $I\!K$. Due to the duality relation between the radial function of a convex body and the support function of its polar body, given by $h_{K^*}=1/ \rho_K$ (see \cite{G}, page 20), our problem is to find and study a measure $\mu$ such that $(\rho_{I\!K})^{-1}=C \mu$. Following the ideas in  \cite{L}, we will use a well-known result of Goodey and Weyl \cite{GW} relating the Radon and the cosine transforms: If $\Delta_n$ is the spherical Laplacian on $S^{n-1}$, then
\begin{equation}
 \label{gw}
     C^{-1}=\frac{1}{\omega_{n-1}} \left(  \Delta_n+n-1 \right) R^{-1},
\end{equation}
where $\omega_{n-1}$ is the Lebesgue measure of the unit sphere in $\mathbb{R}^{n-1}$. In the spherical coordinates $(\sqrt{1-t^2}\,\xi,t)$, $0 \leq t \leq 1$, $\xi \in S^{n-2}$, the spherical Laplacian is expressed by
\[
     \Delta_{n}=(1-t^2) \frac{\partial^2}{\partial t^2}- (n-1) t \frac{\partial}{\partial t}+\frac{1}{1-t^2} \Delta_{n-1},
\]
and $ \Delta_{n-1}$ is applied only to $\xi$, not to $t$ (see \cite{L}, page 9). Hence, for rotationally symmetric functions, equation (\ref{gw}) can be written as 
 \begin{equation}
   \label{laplaceb}
       C^{-1}= \frac{1}{\omega_{n-1}} \left( (1-t^2)\frac{d^2}{dt^2}-(n-1)t\frac{d}{dt}+(n-1)\, I\!d \right) R^{-1}.
 \end{equation}
To simplify the notation, we will write $\Box=\left( (1-t^2)\frac{d^2}{dt^2}-(n-1)t\frac{d}{dt}+(n-1)\, I\!d \right)$.

Now we need the formulas for the spherical Radon transform and its inverse, when acting on rotationally symmetric functions (the derivation of these formulas can be found in the book  \cite{G}, Theorems C.2.9 and C.2.10, p. 432).  Let $f,g$ be rotationally symmetric functions on $S^{n-1}$ such that $f=Rg$. Then 
\begin{equation}
 \label{rad}
   f(\arcsin x)=\frac{C_n}{x^{n-3}} \int_0^x  g(\arccos t) (x^2-t^2)^{(n-4)/2}\, dt,
\end{equation}
for $0<x\leq 1$. The inversion formula is 
\begin{equation}
  \label{invradon2}
   g(\arccos t)=\widetilde{C}_n \, t \left( \frac{1}{t} \frac{d}{dt} \right)^{n-2} 
        \int_0^t  f(\arcsin  x) \, x^{n-2} (t^2-x^2)^{(n-4)/2} \, dx,
\end{equation}
for $0 < t \leq 1$. In all following calculations, we will omit the constants $\omega_{n-1},C_n,  \widetilde{C}_n$ in formulas (\ref{laplaceb}), (\ref{rad}) and (\ref{invradon2}), since their effect is just to dilate the bodies we are considering. Setting $g(\arccos t)=\rho_K^{n-1}(t)$  and $f(\arcsin x)=\rho_{I\!K}(x)$ in (\ref{rad}), we have
\begin{equation} 
  \label{rok}
     \rho_{I\!K}(x)=\frac{1}{x^{n-3}} \int_0^x \rho_K^{n-1}(t) \left( x^2-t^2 \right)^{(n-4)/2}  \, dt,
\end{equation}
for $0 < x \leq 1$. Observe that in (\ref{rok}) the function $\rho_K$ is written in terms of $t$, the cosine of the vertical angle, while $\rho_{I\!K}$ is a function of $x$, the sine of the vertical angle. 

Define 
\begin{equation}
  \label{h}
    h_n(x)=\int_0^x \rho_K^{n-1}(t) \left( x^2-t^2 \right)^{(n-4)/2}  \, dt.
\end{equation}
Then the body $I\!K$ is not a polar zonoid if 
\begin{equation}
 \label{invcos}
     \Box  R^{-1} \left(\frac{x^{n-3}}{h_n(x)} \right)  (t)
\end{equation}
is a  negative measure for some $t \in [0,1]$. 

Let us discuss the regularity of  (\ref{invcos}).  Our starting function $\rho_K$ is continuous and non-zero on $S^{n-1}$, since it is the radial function of a convex body. In  \cite{A}, it was proven that if $\rho_K$ is continuous, then $\rho_{IK}$ is at least in $C^{(n-2)/2}$ at most points (the regularity increase is different at the axis of revolution and at the equator, see the above cited paper for the details). When we apply  the inverse Radon transform to  $\rho_{I\!K}^{-1}$,  the regularity decreases  by the same amount, namely  $(n-2)/2$. As a result, $R^{-1} \left( \frac{x^{n-3}}{h_n(x)} \right) $ is once again a continuous, non-negative function at all points. Next, we apply the operator $\Box$, which  contains a second order derivative. If our initial radial function $\rho_K$ is in the class $C^2(S^{n-1})$,  then (\ref{invcos}) is a continuous function. If $\rho_K$ is $C^1(S^{n-1})$,  then (\ref{invcos}) is a piecewise continuous functions, with jumps at the points where $\rho_K$ is $C^1$ but not $C^2$. Finally,  at any point $t_0$ where $\rho_K$ is continuous but not $C^1$, two differentiations will result in a delta measure supported at the point $t_0$. In all cases, we obtain that (\ref{invcos}) is a measure. 

 We are now ready to study (\ref{invcos}) in the 4-dimensional case (Section 3) and when $n=6$ (Section 4).


\section{Dimension four: A flat-top condition}

As mentioned in the introduction, in dimension four all the operators we are considering have very simple forms. In particular, if $f=Rg$, the inversion formula for the Radon transform  (\ref{invradon2}) simplifies to 
\begin{equation}
  \label{invrad4}
     g(\arccos t)= \frac{d}{dt} (t f(\arcsin t)). 
\end{equation}
Hence, by (\ref{invcos}) we have to compute
\[
    \Box  \left( \frac{d}{dt}  \left( \frac{ t^2}{ h_4(t)} \right)  \right).
\]
In the rest of this section we will write $h(t)$ instead of $h_4(t)$ for notational convenience. Computing the derivative with respect to $t$, we  obtain
\begin{equation}
  \label{defg}
 \Box \left( \frac{2t}{h(t)}-\frac{t^2h'(t)}{h^2(t)} \right)=\Box g(t) = (1-t^2) g''(t) -3t g'(t)+3g(t).
\end{equation}
If we assume that $g$ is of class $C^2$ in a neighborhood of $t=1$ (or equivalently, $\rho_K \in C^2$), and we evaluate the expresion in (\ref{defg}) at  $t=1$, we obtain the following local condition:

\begin{center}{\bf
 If $\mathbf{ g'(1)-g(1)>0}$, then $I\!K$ is not a polar zonoid. }
\end{center}


Computing the derivative of $g$ and simplifying, $g'(1)-g(1)>0$ is equivalent to 
\[
     -3h(1)h'(1)-h(1)h''(1)+2(h'(1))^2>0.
\]
Now recall that $h(x)=\int_0^x \rho_K^3(t) \, dt$. Then $h'(1)=\rho_K^3(1)$, $h''(1)=3\rho_K^2(1)\rho_K'(1)$, and the above condition can be rewritten as
\[
       2\rho_K^4(1)>3h(1) \left( \rho_K(1)+\rho_K'(1) \right).
\]
We have proven:

\begin{Proposition}
 \label{flattop4full}
Let $K$ is a four-dimensional origin-symmetric convex body of revolution whose radial function $\rho_K(t)$ is of class $C^2$ in a neighborhood of the axis of revolution (which corresponds to $t=1$). If the following condition holds,
\begin{equation}
  \label{nec1}
       2\rho_K^4(1)>3 \left( \int_0^1  \rho_K^{3}(t)  \, dt  \right) \left( \rho_K(1)+\rho_K'(1) \right),
\end{equation}
then the intersection body of $K$ is not a polar zonoid. 

\end{Proposition}

Since $\rho_K$ is strictly positive, (\ref{nec1}) will certainly  be satisfied if $\rho_K(1)+\rho_K'(1)=0$. Geometrically, this means that  $K$ is ``flat" at the axis of revolution. For example, the bodies of revolution in $\mathbb{R}^n$ obtained by rotation of the two-dimensional $\ell^p$ balls around the vertical axis satisfy this condition if $p>2$. Also, a body of revolution that has an ($n-1$)-dimensional face perpendicular to the axis of revolution satisfies $\rho_K(1)+\rho_K'(1)=0$.  Hence we have the following criterion:

\begin{Corollary} 
 \label{flattop4}
Let $K$ be a four-dimensional origin-symmetric convex body of revolution, such that $\rho_K(t)$ is of class $C^2$ in a neighborhood of the axis of revolution, and satisfies  $\rho_K(1)+\rho_K'(1)=0$. Then the intersection body of $K$ is not a polar zonoid. 
\end{Corollary}

\begin{figure}[h!]
 \begin{center}
 \includegraphics[width=3.5in]{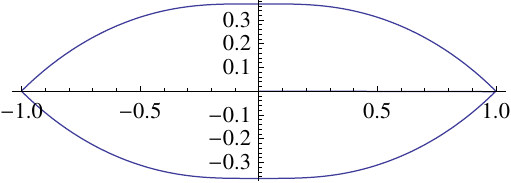}
 \end{center}
\caption{\small A body satisfying the ``flat-top" condition. }
\label{flatt}
\end{figure}

\begin{figure}[h!]
 \begin{center}
 \includegraphics[width=1in]{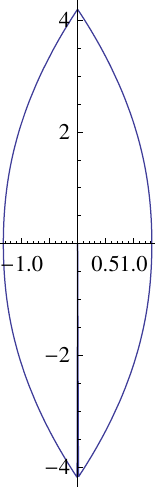}
 \end{center}
\caption{\small The four-dimensional intersection body of the body in Figure \ref{flatt}, which is not a polar zonoid.}
\label{flatt2}
\end{figure}

Flat-top bodies are easy to find, and thus Corollary \ref{flattop4} provides us with many examples of four-dimensional bodies of revolution whose intersection bodies are not polar zonoids. Figure  \ref{flatt} shows the body with radial function $e^{-t}$ for $0 \leq t \leq 1$, which satisfies the flat-top condition, and Figure \ref{flatt2} shows its intersection body in $\mathbb{R}^4$. 

\newpage

The next example is a body that does not have a flat top, but nonetheless verifies the conditions of Proposition \ref{flattop4full}. 

\begin{ex} 
\label{cylcap}
\end{ex}
We consider a four-dimensional cylinder of radius 1/2, with two spherical caps also of radius 1/2 attached to both ends of the cylinder. The radial function of this body is 
\[
       \rho(t)=\left\{    \begin{array}{cc}
                       \frac{1}{2\sqrt{1-t^2}},  &   0\leq t \leq \frac{1}{\sqrt{2}} \\
                        t,  & \frac{1}{\sqrt{2}} \leq t \leq 1. 
 \end{array}
\right. 
\]
Then $\rho(1)+\rho'(1)=2$, $\displaystyle \int_0^1 \rho^3(t) dt =5/16 $, and hence (\ref{nec1}) is satisfied for this body, which implies that its intersection body is not a polar zonoid. More examples are shown in Appendix A.


\section{A geometric  condition in dimension six}

We now consider the six-dimensional case. By (\ref{invcos}), we have to compute 
\[
  \Box  R^{-1} \left(\frac{x^{3}}{h_6(x)} \right) (t),   
\]
where 
\begin{equation} 
 \label{h6}
h_6(x)= \int_0^x \rho_K^{5}(t) (x^2-t^2)  \, dt. 
\end{equation}
Again, for simplicity, we will write $h$ instead of $h_6$ through this Section.  First we compute the inverse Radon transform of $(x^3/h(x))$, using  (\ref{invradon2}):
\[
   R^{-1}\left(\frac{x^3}{h(x)}\right)(t)=   t \left( \frac{1}{t} \frac{d}{dt} \right)^{4} 
        \int_0^t   \frac{x^{7}}{h(x)} \, (t^2-x^2) \, dx  =2 \, t  \left( \frac{1}{t} \frac{d}{dt} \right)^{3} 
        \int_0^t   \frac{x^{7}}{h(x)} \, dx
\]
\[
     = 2 \, t  \left( \frac{1}{t} \frac{d}{dt} \right)^{2} 
          \frac{t^{6}}{h(t)} \, dx
     = 2 \, \frac{d}{dt}  \left( \frac{6t^4}{h(t)} -\frac{ t^5 h'(t)}{h^2(t)} \right)
\]
\[
    =2 \, \left( \frac{24t^3}{h(t)}-\frac{t^5h''(t)+11t^4h'(t)}{h^2(t)}+\frac{2t^5(h'(t))^2}{h^3(t)}\right)=2 \, g(t).
\]
Ignoring the factor of 2 in the above expression, we now have to calculate
\[
    \Box g(t)= (1-t^2) g''(t) -5t \, g'(t)+5g(t).
\]
As in the four dimensional case, evaluating of the above expression at $t=1$ provides us with a local condition: If $g\in C^2$ in a neighborhood of $t=1$, and $g'(1)-g(1)>0$, then the intersection body of $K$ is not a polar zonoid. Differentiating $g$, we obtain that $g'(1)-g(1)$ is equal to
\begin{equation}
   \label{vamos}
      \frac{1}{h^4(1)} \left( 
         48h^3(1)-6(h'(1))^3+2h(1)h'(1) \left( 15h'(1)+3h''(1) \right) \right.
\end{equation}
\[ 
    \left. -h^2(1) \left( 57h'(1)+15h''(1)+h'''(1) \right)
     \right).
\]
The term $h^4(1)$ in the denominator is positive, so we need the numerator to be positive. From the definition of $h$ (\ref{h6}), and writing $r(t)=\rho_K^5(t)$,  we have that $h'(1)=2 \int_0^1 r(t) dt$,  $ \; h''(1)=h'(1)+2 r(1)$, $\; h'''(1)=4r(1)+2r'(1)$. Thus the numerator in (\ref{vamos}) will be positive if and only if
\[
  h^2(1)\left(5r(1)+r'(1) \right)+ 24\left(\int_0^1 t^2 r(t) \,dt \right)^3 < 12h(1)r(1)\left(\int_0^1 t^2 r(t) \,dt \right).
\]
We summarize the result in the following proposition.

\begin{Proposition}
 \label{nec16}
Let $K$ be a six-dimensional origin-symmetric convex body of revolution with radial function $\rho
_K(t) $ of class $C^2$ in a neighborhood of  $t=1$.  Let $h(x)= \int_0^x \rho_K^{5}(t) (x^2-t^2)  \, dt$,  $\displaystyle k(1)=\int_0^1 t^2 \rho_K^5(t) \,dt$, and  $r(t)=\rho_K^5(t)$. If 
\begin{equation}
 \label{nc6}
   h^2(1)\left(5r(1)+r'(1) \right)+ 24k^3(1) < 12h(1)k(1)r(1),
\end{equation}
then $I\!K$ is not a polar zonoid.
\end{Proposition}

Note that $5r(1)+r'(1)=0$ is equivalent to $\rho_K(1)+\rho_K'(1)=0$. Thus we have the following condition for bodies with a flat top.

\begin{Corollary} 
 \label{flattop6}
Let $K$ be a six-dimensional origin-symmetric convex body of revolution such that $\rho_K(t)$ is of class $C^2$ in a neighborhood of $t=1$, and $\rho_K(1)+\rho_K'(1)=0$. Let $h,r,k$ be defined as in Proposition \ref{nec16}. If \begin{equation}
 \label{nc6flat}
        2k^2(1)<h(1)r(1),
\end{equation}
then the intersection body of $K$ is not a polar zonoid.
\end{Corollary}

It is interesting to note the difference between Corollary \ref{flattop4} and Corollary \ref{flattop6}.  In dimension four, a flat top on $K$ is enough to guarantee that $I\!K$ is not a polar zonoid. In dimension six, a flat top is no longer enough, and we are required to add condition (\ref{nc6flat}). Intuitively, this condition means that the body $K$ is ``fat" around the 
 and ``thinner" close to the pole. This is a result of the presence of the weights $(1-t^2)$ and $t^2$  in the definitions of $h(1)$ and $k(1)$. Example \ref{3bodies}   illustrates this.

\begin{figure}[h!]
 \begin{center}
 \includegraphics[width=2.1in]{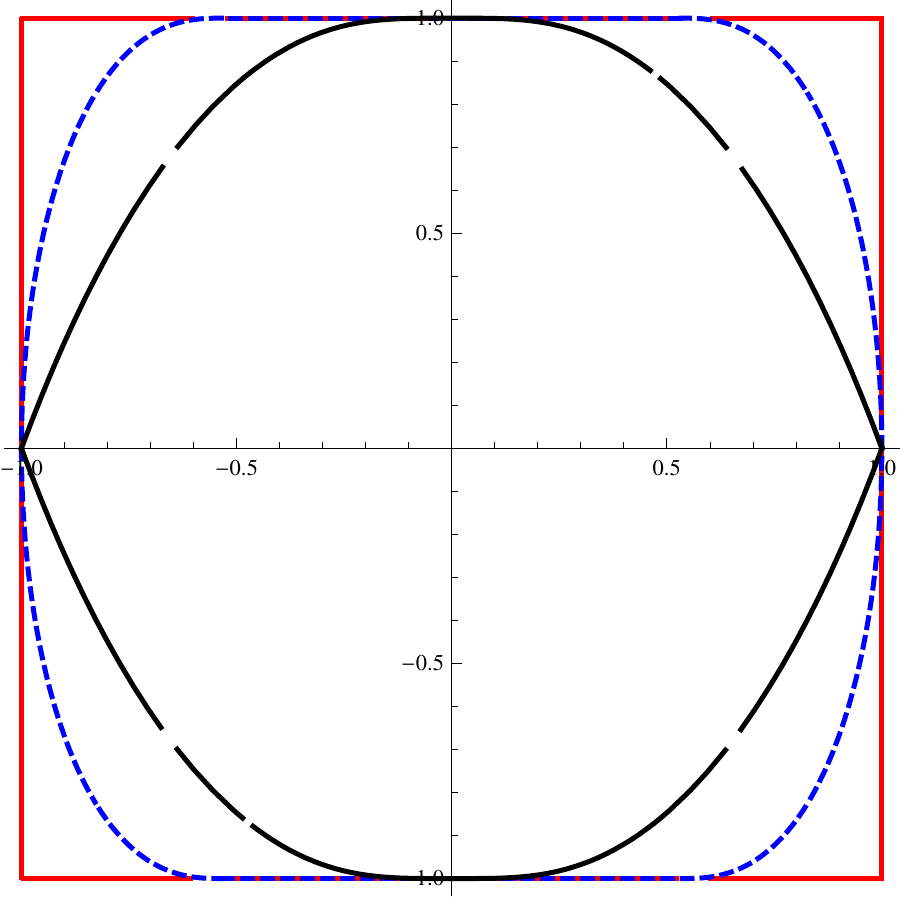}
 \end{center}
\caption{\small Comparison of three bodies.}
\label{comparison}
\end{figure}

\begin{ex}
 \label{3bodies}
\end{ex}
Figure \ref{comparison} shows the cross-section of three bodies of revolution in $\mathbb{R}^6$: The cylinder $C$, whose radial function is 
\[
   \rho_C(t)=\left\{ 
 \begin{array}{cc}
     \frac{1}{\sqrt{1-t^2}} &   0\leq t \leq \sqrt{2}/2 \\
     \frac{1}{t}  &\sqrt{2}/2 \leq t \leq 1
 \end{array}
 \right.;
\]
the dashed body $L$, which is an affine transformation of the body with radial function 
\[
   \rho_L(t)=\left\{ 
 \begin{array}{cc}
   \frac{3 - 16 (1 - t^2) + 28(1 - t^2)^2}{8(1 - t^2)^{5/2}} &  0\leq t \leq \sqrt{2}/2 \\
    1/t  &\sqrt{2}/2 \leq t \leq 1
 \end{array}
 \right.;
\]
and the smaller body $E$, which is an affine transformation of the body with radial function $\rho_E(t)=e^{-t}$. It is immediate to check that these three bodies satisfy the flat top condition.

 We will now study condition (\ref{nc6flat}) for each of them. The computations are done with {\it Mathematica}.

For the cylinder $C$, we have that $r(1)=1$, 
\[
    h(1)=\int_0^{1/\sqrt{2}} (1-t^2)^{-3/2} \, dt + \int_{1/\sqrt{2}}^1 \frac{1-t^2}{t^5} \, dt=\frac{1}{2} + \frac{3 \pi}{32},
\]
and
\[
    k(1)=\int_0^{1/\sqrt{2}} \frac{t^2}{(1-t^2)^{5/2}} \, dt + \int_{1/\sqrt{2}}^1 t^{-3} \, dt=\frac{5}{6}.
\]
Since $2(5/6)^2>\frac{1}{2} + \frac{3 \pi}{32}$, the cylinder does not satisfy (\ref{nc6flat}).\\

Similarly, for $L$ we obtain $\displaystyle r(1)=1$, $ h(1)=\frac{44239925}{3879876}$, and $ k(1)=\frac{30712575}{14872858}$. In this case, $\displaystyle r(1)h(1)>2(k(1))^2$. By Corollary \ref{flattop6}, $I\!L$ is not a polar zonoid.  Finally, for $E$ we have $\displaystyle r(1)=e^{-1}$, $ h(1)=\frac{23 + 12e^{-5}}{125}$, $ k(1)=\frac{2 -37e^{-5}}{125}$. Again,  $r(1)h(1)>2(k(1))^2$ and thus $I\!E$ is not a polar zonoid.

\begin{remark}
 \label{pert}
\end{remark} From the definitions of $h(1)$ and $k(1)$, the following result is immediate: If $\rho_K$ satisfies (\ref{nc6flat}) and $L$ is a body such that $\rho_L(t) \leq \rho_K(t)$ for $0 \leq t \leq 1/\sqrt{2}$, and $\rho_L(t) \geq \rho_K(t)$ for $1/\sqrt{2} \leq t \leq 1$, then $L$ will also satisfy (\ref{nc6flat}). It is thus quite straightforward to construct families of bodies satisfying the conditions of Corollary \ref{flattop6}. More examples are shown in Appendix B. 

\begin{ex}
\end{ex}
As a further application of Corollary \ref{flattop6}, we study the bodies $K^6_p \subset \mathbb{R}^6$, obtained by rotation a unit ball of two-dimensional $\ell^p$ about the vertical axis. As observed before Corollary (\ref{flattop4}), if $p>2$ then $K_p^6$ verifies the flat top condition $\rho_{K_p^6}(1)+\rho_{K_p^6}'(1)=0$. We checked with the help of the computer that  $K_p^6$ satisfies (\ref{nc6flat}) for $2<p \leq 9.5$. Hence, for $p$ in this range,   $I\!K_p^6$ is not a polar zonoid by Corollary \ref{flattop6}. For $p\geq 9.6$, condition  (\ref{nc6flat}) fails and the Corollary does not allow us to conclude anything. We have no geometric explanation for the change between $9.5$ and $9.6$. 

\bigskip

As the previous example shows, the applications of Corollary  \ref{flattop6} are limited. Although the condition (\ref{nc6flat}) does not for the cylinder (which is the case  $K=K^6_\infty$  in Example 8), we have computed the function $\Box R^{-1} (\rho_K^{-1})(x)$ for all values of $x \in [0,1]$, and found it to be sign changing (and, of course, positive at $x=1$, since condition (\ref{nc6flat}) fails). The complete calculations are included in Appendix C.  Thus, the intersection body of the six-dimensional cylinder is not a polar zonoid, but we cannot obtain this information just from Corollary \ref{flattop6}. On the other hand, as observed in Remark \ref{pert}, once we find a body that satisfies the conditions of the Corollary, we can perturb it to construct many other examples. \\

The approach taken in this paper provides conditions that become harder to compute as the dimension increases. We have calculated $g'(1)-g(1)$  for $n=8$, but we were not able to find its geometric meaning, even with the extra assumption $\rho_K(1)+\rho_K'(1)=0$. Nonetheless, we feel that at least in dimension six the approach was useful in providing a characterization and many new examples of intersection bodies that are not polar zonoids. 


In the three appendices, we present several applications of Proposition \ref{flattop4full} and Corollary \ref{flattop6}. The computations in these sections have been done with {\it Mathematica}.

\bigskip

\section*{Appendix A: Application of Proposition  \ref{flattop4full} to a family of cylinders with spherical caps in $\mathbb{R}^4$.}

\begin{figure}[h!]
 \begin{center}
 \includegraphics[width=1.2in]{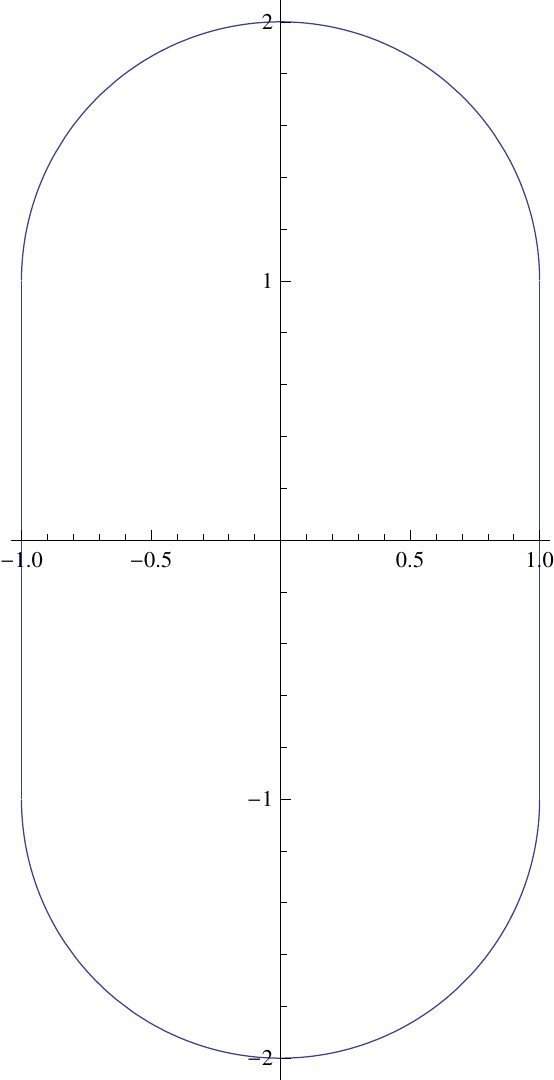}
 \end{center}
\caption{\small The body $K_1$ in Appendix A}
\label{k1}
\end{figure}

\begin{figure}[h!]
 \begin{center}
 \includegraphics[width=1.6in]{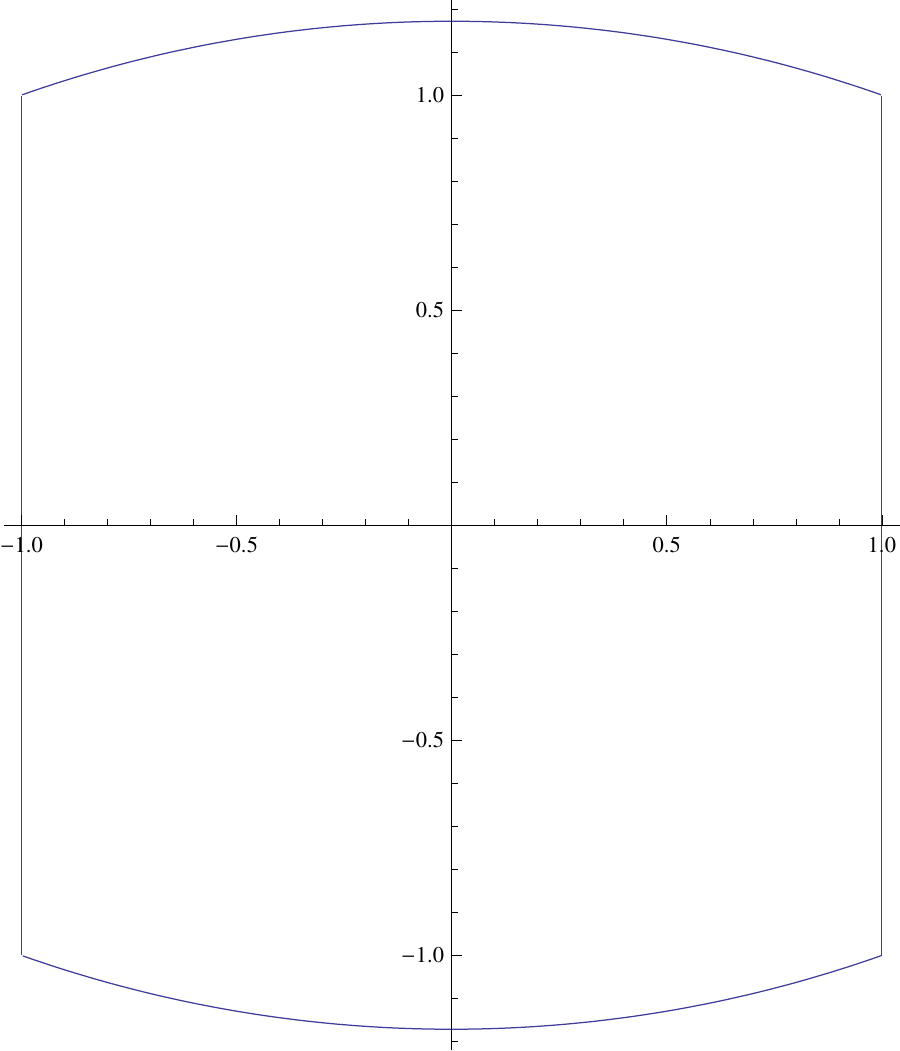}
 \end{center}
\caption{\small The body $K_3$ in Appendix A}
\label{k3}
\end{figure}

Inspired by Example \ref{cylcap} in Section 3, in Appendix A we check the hypotheses of Proposition \ref{flattop4full} for a family of origin-symmetric  bodies of revolution $K_M$, each of them being a cylinder centered at the origin and with height 2, and having two spherical caps attached. In two dimensions, the top cap is an arc of the circle with center at the point $(0,1-\sqrt{M^2-1})$ and radius $M \geq 1$. When $M=1$, $K_1$ is a dilate of the body in Example \ref{cylcap}. As $M$  tends to infinity, the top becomes flatter and $K_M$ tends to the cylinder in the radial metric. Their radial function  is given by 
\[
       \rho_M(t)=\left\{    \begin{array}{cc}
                      (1-t^2)^{-1/2},  &   0\leq t \leq \frac{1}{\sqrt{2}} \\
                      t (1 - \sqrt{M^2 - 1}) + 
  \sqrt{t^2(M^2 - 2 \sqrt{M^2 - 1}) + 2 \sqrt{M^2 - 1}},  & \frac{1}{\sqrt{2}} \leq t \leq 1.
 \end{array}
\right. 
\]
 Two of these bodies are shown in Figures \ref{k1}, \ref{k3}.

Observe that the body $M$ has boundary at least of class $C^2$ at every point except at $t=1/\sqrt{2}$, where it is only continuous. The intersection body $I\!K_M$ looks like a barrel. The flat top condition in Corollary \ref{flattop4} does not hold for any finite $M$. However, we have seen in Example \ref{cylcap} that $K_1$ satisfies the  condition (\ref{nec1}) in Proposition \ref{flattop4full}, and our current goal is to check it for other values of $M$. We thus define the function $w(M)= 2\rho_{K_M}^4(1)-3 \left( \int_0^1  \rho_{K_M}^{3}(t)  \, dt  \right) \left( \rho_{K_M}(1)+\rho_{K_M}'(1) \right)$, and we need to check for which values of $M$ we have $w(M)>0$.

With {\it Mathematica}, we compute  $\rho_M(1)=1+M-\sqrt{M^2-1}$, 
\[
    \rho_M(1)+\rho'_M(1)=\frac{2 (1 +M) (M - \sqrt{ M^2-1})}{M}, 
\]
and $ \int_0^1 \rho_M^3(t)\, dt=1 - \frac{1}{(4 (M^2-1)^{3/4} (M^2 - 2 \sqrt{M^2-1})^2)} \left[ 3 M^8 (M^2-1)^{1/4} - 16 M^3 (M^2-1)^{3/4} \right.$

\noindent $ - 
   8 (M^2-1)^{1/4} \sqrt{M^2 + 2 \sqrt{M^2-1}}  +   4 M^7 ((M^2-1 )^{3/4}-4 (M^2-1)^{1/4}) +  16 M^5 ((M^2-1)^{1/4} $

\noindent $ + (M^2-1)^{3/4})  +   M^4 (36 ( M^2-1)^{1/4}
-  6 (M^2-1)^{1/4} \sqrt{M^2 + 2 \sqrt{M^2-1}}$

\noindent $ +      8 ( M^2-1)^{3/4} \sqrt{M^2 + 2 \sqrt{M^2-1}} - 4 \sqrt{-2 + M^2 (2 + \sqrt{M^2-1})}$

\noindent$-     4 \sqrt{(M^2-1) (-2 + M^2 (2 + \sqrt{M^2-1}))})
-   4 M^2 (3 (M^2-1)^{1/4} + 3 (M^2-1)^{3/4}$

\noindent$ -  4 (M^2-1)^{1/4}\sqrt{M^2 + 2 \sqrt{M^2-1}}
 + 2 (M^2-1)^{3/4} \sqrt{M^2+ 2 \sqrt{M^2-1}}$

\noindent $ -\sqrt{(M^2-1) (-2 + M^2 (2 + \sqrt{M^2-1}))}) $

\noindent$+    M^6 \left\{ -21 (M^2-1)^{1/4} + 15 (M^2-1)^{3/4} \right. -2 (M^2-1)^{1/4} \sqrt{M^2 + 2 \sqrt{M^2-1}}$

$+ \left.
      4 \sqrt{-2 + M^2 (2 + \sqrt{M^2-1})}  \left. -\sqrt{( M^2-1) (-2 + M^2 (2 + \sqrt{M^2-1}))} \right\}\right]$.

\bigskip

\begin{figure}[h!]
 \begin{center}
 \includegraphics[width=3.5in]{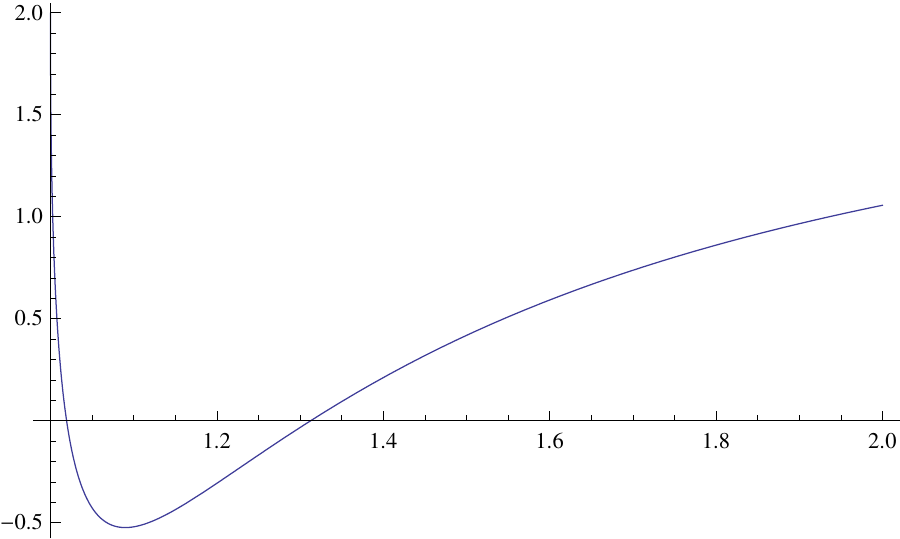}
 \end{center}
\caption{\small The function $w(M)$ for the family of bodies $K_M$. }
\label{fm}
\end{figure}

 Figure \ref{fm} shows the graph of $w$. It is a continuous function,  with exactly two zeros $M_1 \in(1.01942,1.01943)$ and $M_2 \in (1.31290,1.31291)$. The limit of $w(M)$ as $M$ tends to infinity equals 2. By Proposition \ref{flattop4full}, for all the values of $M$ such that $w(M)>0$, we know that $I\!K_M$ is not a polar zonoid. 

\section*{Appendix B: The conditions of Corollary  \ref{flattop6} for a family of bodies of revolution in $\mathbb{R}^6$.}

\begin{figure}[h!]
 \begin{center}
 \includegraphics[width=2.4in]{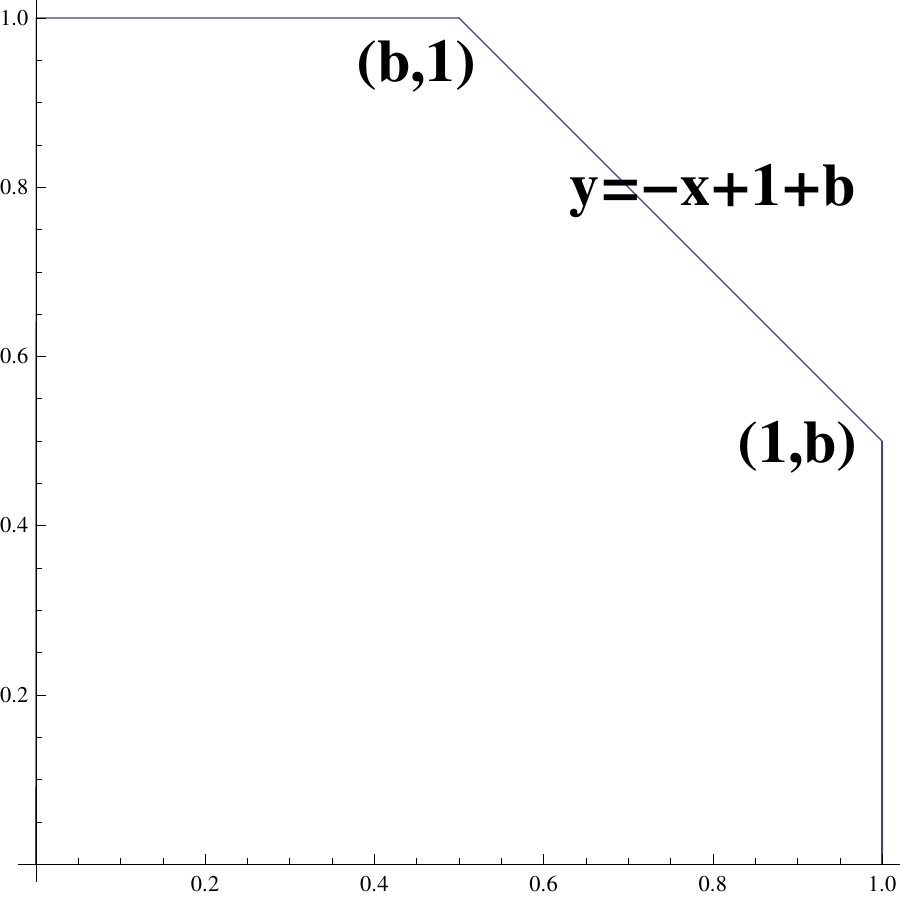}
 \end{center}
\caption{\small $\rho_{K_b}(t)$, $0 \leq t \leq 1$, with $b=1/2$.}
\label{faml}
\end{figure}

Here we will study a family of bodies of revolution in $\mathbb{R}^6$ satisfying the hypotheses of Corollary \ref{flattop6}. For $b \in [0,1]$,  $K_b$ is the body of revolution of an octagon whose sides depend on the parameter $b$. Figure \ref{faml} shows the part of the octagon in the first quadrant.  For $b=0,1$, the octagon becomes a square, and $K_0$, $K_1$ are, respectively, a double cone and a cylinder. The radial function of $K_b$ is given by
\[ 
  \rho_{K_b}(t)=\left\{    \begin{array}{cc}
                      1/\sqrt{1-t^2} ,  &   0\leq t \leq \frac{b}{\sqrt{1+b^2}} \\
                    \frac{1+b}{t+\sqrt{1-t^2}} &   \frac{b}{\sqrt{1+b^2}} \leq t \leq  \frac{1}{\sqrt{1+b^2}}\\
                     1/t,  & \frac{1}{\sqrt{1+b^2}} \leq t \leq 1. 
 \end{array}
\right. 
\]

\begin{figure}[h!]
 \begin{center}
 \includegraphics[width=3in]{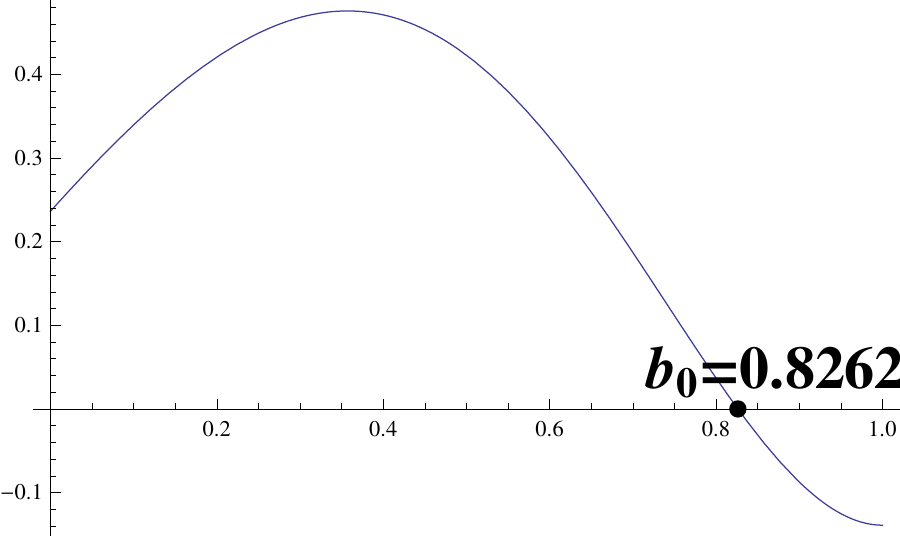}
 \end{center}
\caption{\small Condition (\ref{nc6flat}) for the bodies $K_b$ is satisfied if $b \in [0,b_0]$.}
\label{pol}
\end{figure}

For $b>0$, the bodies $K_b$ are $C^2$ in a neighborhood of $t=1$ and satisfy the flat top condition $\rho_{K_b}(1)+\rho_{K_b}'(1)=0$. We will now check  condition (\ref{nc6flat}). The condition will be satisfied if $h(1)-2(k(1))^2>0$, where
\[
     h(1)=\int_0^1 (1-t^2) \rho_{K_b}^5(t) \, dt = \frac{1 + 5 b - b^5}{4},
\]
\[
     k(1)=\int_0^1 t^2 \rho_{K_b}^5(t) \, dt =\frac{1 + 5 b + 10 b^2 - 5 b^4 - b^5}{12}. 
\]
Thus, $h(1)-2(k(1))^2$ is a tenth-degree polynomial in the variable $b$. Its graph is shown in Figure \ref{pol}. It has only one root $b_0=0.826279...$ in the interval $[0,1]$. By Corollary \ref{flattop6}, for $0\leq b <b_0$, the body $I\!K_b$ is not a polar zonoid.

\section*{Appendix C: The intersection body of the cylinder in $\mathbb{R}^6$.}

\begin{figure}[h!]
 \begin{center}
 \includegraphics[width=3in]{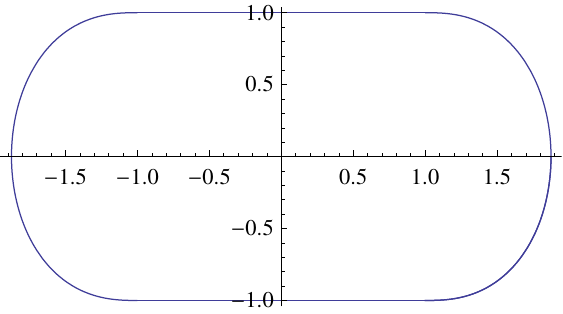}
 \end{center}
\caption{\small The intersection body of the cylinder in $\mathbb{R}^6$.}
\label{ibc}
\end{figure}

The centered cylinder $C$ with radius 1 and height 2 has radial function given, in terms of $t$ (the cosine of the vertical angle) by
\[ 
  \rho_C(t)=\left\{    \begin{array}{cc}
                      1/\sqrt{1-t^2} ,  &   0\leq t \leq \frac{1}{\sqrt{2}} \\
                     1/t,  & \frac{1}{\sqrt{2}} \leq t \leq 1. 
 \end{array}
\right. 
\]
From (\ref{rok}), the intersection body of $C$ in $\mathbb{R}^6$ has radial function (in terms of $x$, the sine of the vertical angle), given by 
\[
    \rho_{I\!C}(x)=\left\{    \begin{array}{cc}
             1/\sqrt{1-x^2} ,  &   0\leq x \leq \frac{1}{\sqrt{2}} \\
                (3 - 16 x^2 + 28 x^4)/(8x^5),     & \frac{1}{\sqrt{2}} \leq x \leq 1. 
 \end{array}
\right. 
\]
See Figure \ref{ibc}.  In this Appendix we compute $\Box \left( R^{-1} (\rho_{IC}^{-1}) (t) \right)$ for all values of $t$.

\begin{figure}[h!]
 \begin{center}
 \includegraphics[width=3.2in]{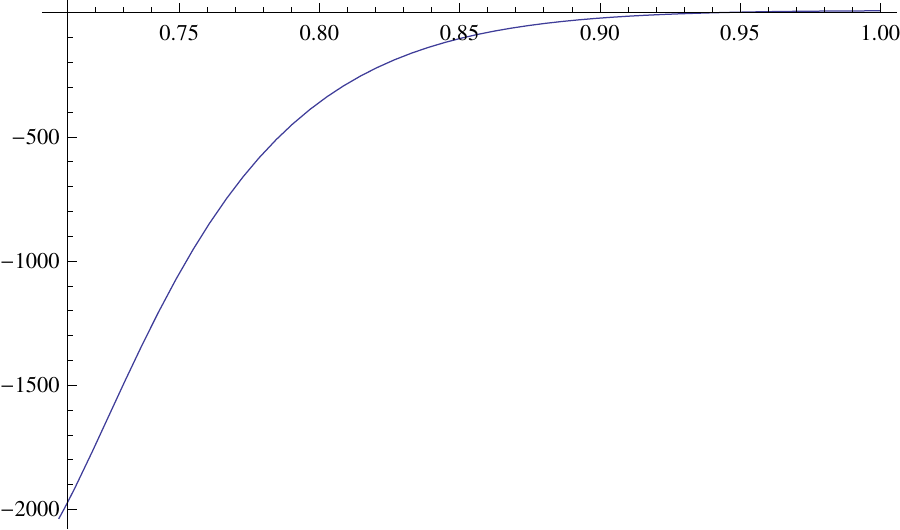}
 \end{center}
\caption{\small The function  $\Box \left( R^{-1} (\rho_{IC}^{-1}) \right)$ for the cylinder}
\label{cyl6}
\end{figure}

 First we invert the Radon transform using (\ref{invradon2}), and we obtain the function $g(t)=R^{-1} (\rho^{-1}_{I\!C})(t)$
\[
    = \left\{    \begin{array}{cc}
            (6 - 24 t^2 + 16 t^4)/(1 - t^2)^{3/2} ,  &   0 \leq t \leq \frac{1}{\sqrt{2}} \\
                (256 t^5 (27 - 192 t^2 + 510 t^4 - 672 t^6 + 392 t^8))/(3 - 16 t^2 + 
  28 t^4)^3  ,   & \frac{1}{\sqrt{2}}  \leq t  \leq 1. 
 \end{array}
\right. 
\]
Observe that $g$ is continuous but not differentiable at $t=1/\sqrt{2}$. Its first and second derivatives are given by 
\[
    g'(t)=\left\{    \begin{array}{cc}
          -(2 t (15 - 20 t^2 + 8 t^4))/(1 - t^2)^{5/2} ,  &   0 \leq t < \frac{1}{\sqrt{2}} \\
               (256 t^4 f_1(t))/(3 - 16 t^2 + 28 t^4)^4 ,  & \frac{1}{\sqrt{2}}  < t  \leq 1. 
 \end{array}
\right. 
\]
where $f_1(t)=405 - 3600 t^2 + 11550 t^4 - 19776 t^6 + 26208 t^8 - 
   25088 t^{10} + 10976 t^{12}$, and 
\[
      g''(t)=c\, \delta_{1/\sqrt{2}}+\left\{    \begin{array}{cc}
                    -(30/(1 - t^2)^{7/2}),   &   0 \leq t < \frac{1}{\sqrt{2}} \\
     (15360 t^3 f_2(t))/(3 - 16 t^2 + 28 t^4)^5, & \frac{1}{\sqrt{2}}  < t  \leq 1. 
 \end{array}
\right. 
\]
where $f_2(t)=81 - 648 t^2 + 432 t^4 + 6912 t^6 - 16848 t^8 + 9856 t^{10}$. The value of the constant $c$ is 
\[
 c=\lim_{t \rightarrow 1/\sqrt{2}^+} g'(t)- \lim_{t \rightarrow 1/\sqrt{2}^-} g'(t)=184-(-56)=240.
\]
Finally, we put everything together. The function  $(1-t^2)g''(t)-5tg'(t)+5g(t)$ is identically equal to zero for $0\leq t < 1/\sqrt{2}$. For $1/\sqrt{2}  < t  \leq 1$, it is equal to 
\[
  (46080 t^3 (27 - 270 t^2 + 720 t^4 + 240 t^6 - 2800 t^8 + 
   2208 t^{10}))/(3 - 16 t^2 + 28 t^4)^5.
\] 
The graph of this function is shown in Figure \ref{cyl6}. It takes negative values, and thus $I\!C$ is not a polar zonoid. At the point $t=1$  the function is positive, as we already knew because the condition (\ref{nc6flat}) fails for the cylinder.

{\it Acknowledgements:} The author wishes to thank D. Ryabogin and A. Koldobsky for useful  conversations about this work, and also the referee for his suggestions to improve the paper.

\end{document}